\documentclass[12pt,twoside]{article}

\usepackage{amssymb}
\usepackage{pst-node,pstcol,pstricks,multido,pst-plot,pst-text,pst-3d}%
\usepackage{graphicx}
\usepackage{amsfonts}
\newtheorem{example}{Example}[section]

\newtheorem{theorem}[example]{Theorem}

\newtheorem{corollary}[example]{Corollary}

\newtheorem{proposition}[example]{Proposition}

\def\Prod{\rm prod}

\def\<{\langle}
\def\>{\rangle}

\def\N{{\mathbb N}}
\def\M{{\mathbb M}}
\def\F{{\mathbb F}}

\def\H{{\sf H}}

\def\Card{{\rm Card\ }}
\def\cqfd{$\ \ \ \blacksquare$\\ \\}
\def\Alph{\mbox{\rm Alph}}

 \font\tengoth=eufb10
\font\sevengoth=eufb7 \font\fivegoth=eufb5

\newfam\gothfam
\textfont\gothfam=\tengoth \scriptfont\gothfam=\sevengoth
\scriptscriptfont\gothfam=\fivegoth
\def\goth{\fam\gothfam\tengoth}

\def\ashuff#1#2#3{
\kern 1pt \vrule height#1 \overline{\vrule height#3 width 0pt
\hskip#2} \rule{.3pt}{#1}\overline{\vrule height#3 width 0pt
\hskip#2} \rule{.3pt}{#1} \kern 1pt }

\title{Transitive Hall sets}

\author{{\bf 
G.H.E. Duchamp$^\dag$, M. Deboysson-Flouret$^\ddag$, J.-G. 
Luque$^\star$}\\ 		
($^\dag$) LIPN, UMR CNRS 7030\\
Institut Galil\'ee - Universit\'e Paris-Nord\\ 99, avenue
Jean-Baptiste Cl\'ement\\ 93430 Villetaneuse, France\\
{\normalsize \tt ghed@lipn.univ-paris13.fr}\\
$(^\ddag)$ Laboratoire d'Informatique du Havre (EA3219)\\
25 rue Philippe Lebon,  BP 540\\ 
76058 Le Havre cedex\\
{\normalsize \tt marianne.deboysson@univ-lehavre.fr }\\
$(^\star)$ IGM,   Laboratoire d'informatique UMR 8049 IGM-LabInfo\\ 77454
Marne-la-Vall\'ee Cedex 2, France.\\ {\normalsize \tt
 luque@univ-mlv.fr} \\
}

\date{}

\begin{document}
 \maketitle

\begin{center}Dedicated to G\'erard X. Viennot\end{center}
 \begin{abstract}
We give the definition of Lazard and Hall sets in the context of
transitive factorizations of free monoids. The equivalence of the
two properties is proved. This allows to build new effective bases
of free partially commutative Lie algebras. The commutation graphs
for which such sets exist are completely characterized and we
explicit, in this context, the classical PBW rewriting process.
\end{abstract}

\section{Introduction}
The correspondence with Lie algebras  is not so clear  for monoids
as it is  in the case of groups (Lie groups and, after the work of
Magnus \cite{Cha} and Lazard \cite{Laz} combinatorial groups). In
fact, the first connection between the two realms was done by
means of the notion of factorization \cite{Vi2}.

An ordered family $(M_i)_{i\in I}$ of submonoid is said to be a
factorization if the product mapping $\coprod_{i\in
I}M_i\rightarrow M$ is one to one (see paragraph \ref{schutz}).
For $|I|=2$, one gets the notion of a bisection related to the
flip-flop
% {\tt v\'erifier dans Mathsciencenet s'il a publie en anglais}
\footnote{In french "bascule"} Lie-algebra by Viennot \cite{Vi2}.
On the other end, when all the $M_i$ have a single generator, one
obtains a complete factorisation whose generating series is equal
to the Hilbert series of the Free Lie algebra. It has been shown
in \cite{DK1} that this property stills holds in case of partial
commutations and the link between the Lyndon basis and Lazard
elimination in this context has been elucidated \cite{KL}. Here,
we consider a more general construction: the transitive Hall sets.
\smallskip\\ \\
The structure of the paper is the following. In section
(\ref{fact}), we give a partially commutative version of
Sch\"utzenberger's factorization theorem. We introduce in Section
(\ref{TLS}) the notions of transitive Lazard and transitive Hall
sets in the free algebras of trees. Transitive Lazard sets are
defined by means of iterations of transitive bisections \cite{DL}.
Transitive Hall sets classically use a description of the total
order of the factorization which, here, must be compatible with
the commutations.  We prove that the two notions coincide and
allow the construction of bases of $L(A,\theta)$ as well as a
natural Poincar\'e-Birkhoff-Witt's rewriting process (PBW).

\section{Factorizations of a trace monoid\label{fact}}
Let $A$ be an alphabet and $\theta\in A\times A-\{(a,a)|a\in A\}$
be a symmetric relation. We will denote by $\M(A,\theta)$ a trace
monoid (or a free partially commutative monoid) over the alphabet
$A$ and whose commutations are defined by $ab=ba$ when
$(a,b)\in\theta$.\\ The only way to get the equality
\begin{equation}
\underline{\M(A,\theta)}=\underline{a^*}.\underline{\M(X,\theta_X)}
\end{equation}
where  $X$ a subset of $\M(A,\theta)$ is to impose for each pair
of letters $(z_1,z_2)\in\theta\cap(A-a)^2$ the commutation
$(z_1,a)\in\theta$ or $(z_2,a)\in\theta$. It is a direct
consequence of the transitive factorization theorem (see Section
(\ref{TypeH}) \cite{DL}). In this paper, we use this property to
define transitive Lazard sets for a certain family of trace monoids. Let
us first show some properties about  conjugacy.

\subsection{Roots of conjugacy classes}
Let $m\in \M(A,\theta)$. It is shown  in \cite{DK3} (see also
\cite{Dub}) that the equation $u=t^p$ ($p\geq 1$) has at most one
solution. When it exists, this solution will be denoted by
$^p\sqrt u$. In the same paper \cite{DK3}, it is shown that if
$g$, $r$, and $t$ are tree traces such that $g=r^q=t^p$ with $q$,
$p\in\N$ then it exists a trace $g'$ and $m\in {\rm lcm}(p,q)\N$
such that  $g=g'^m$. Hence, one defines the {\bf  root} $\sqrt g$
of a trace $g$ as the smallest trace $g'$ satisfying $g=g'^m$ (the
integer $m$ will be called the {\bf exponent} of $g$, we denote
$m={\rm ex}(g)$).
\\
\begin{example}
Let us consider the commutation graph
\[(A,\theta)=a-c-b\]
we have $\sqrt{babcac}=bac$ and ${\rm ex}(babcac)=2$.
\end{example}
We recall here the definition of conjugacy due to Duboc and
Choffrut  \cite{CH}. Two traces $t$ and $t'$ are said to be 
conjugate if
there exists a trace $u$ such that $tu=ut'$. Conjugacy is an
equivalence relation which, in turn, is no more that the
restriction to $\M(A,\theta)$ of the conjugacy relation of the
group $\F(A,\theta)$. Exponents and roots are invariant under
conjugacy in the following sense.

\begin{proposition}\label{conjug}
Let $C$ be a conjugacy class,  $f\in C$, $g\in\M(A,\theta)$ and
$p\in\N$ be such that $f=g^p$. For each $f'\in C$, it exists
$g'\in\M(A,\theta)$ such that $f'=g'^p$. Furthermore $g$ and $g'$
are conjugate.
\end{proposition}

As a consequence of Proposition (\ref{conjug}) root $\sqrt C$ of a
conjugacy class $C$ is uniquely defined as the class $\sqrt
C=\{\sqrt g\}_{g\in C}$.

\subsection{Sch\"utzenberger's factorization theorem for traces}
\label{schutz}

Let $(M_i)_{i\in J}$ be a family of monoids. The restricted
product $\coprod_{i\in J} M_i$ is the submonoid of $\prod_{i\in J}
M_i$ of the families with finite support (all but a finite number
of indices are $1_{M_i}$). Now, if the $(M_i)_{i\in J}$ are
submonoids of a given monoid $M$ and if $J$ is totally ordered,
the product $\Prod : M\times M\rightarrow M$ extends to an arrow
$\Prod : \coprod_{i\in J} M_i\rightarrow M$ by means of the
ordered product.\\ A factorization of a monoid $\M$ is an ordered
family of submonoids $(M_i)_{i\in J}$ such that $\Prod :
\coprod_{i\in J} M_i\rightarrow M$ is one to one. It is classical
that a submonoid $\M$ of $\M(A,\theta)$ has a unique minimal
generating set. Hence, a factorization can be characterized by
the family of the generating sets of its components and will be
denoted by $\F=(Y_i)_{i\in I}$ instead of $\F=(\M_i)_{i\in I}$ if
$Y_i=\M_i-\M_i^2$.
\\ \\
The existence of the root of conjugacy classes allows us to extend
Sch\"utzen\-berger's factorization theorem \cite{onfact} to trace
monoids.
\begin{theorem}
Let $\F=(Y_i)_{i\in J}$ be an ordered family of non-commutative
subsets (i.e. for each i and each pair $(x,y)\in Y_i^2$, $x\neq y$
implies $xy\neq yx$) of $\M(A,\theta)$ and $\langle Y_i\rangle$
the submonoid generated by $Y_i$.

\noindent We consider the following assertions :
\begin{enumerate}
\item\label{it1} The mapping {\Prod}  is into.
\item\label{it2} The mapping {\Prod  } is onto.
\item\label{it3} Each monoid  $\langle Y_i\rangle$ is free.
For each conjugacy class $C$ in $\M(A,\theta)$, if $C$ is
connected   ({\it i.e.} if the restriction of the non-commutation
graph to the alphabet $\Alph(t)=\{a\in A|t=uav\}$ is a connected
graph) then it exists an unique $i\in J$ such that $C\cap\langle
Y_i\rangle\neq\emptyset$ and in this case $C\cap\langle
Y_i\rangle$ is a conjugacy class of
 $\langle Y_i\rangle$. If $C$
is not connected, for each $i\in J$, $C\cap\langle
Y_i\rangle=\emptyset$.
\end{enumerate}
Two of the previous assertions imply the third.
\end{theorem}
{\bf Proof} We prove \ref{it1}) and \ref{it2}) imply \ref{it3}) by
means of the examination of the series
\begin{equation}\label{log}
\log\underline{\M(A,\theta)}-\sum\log\underline{\langle
Y_i\rangle}.
\end{equation}
Remarking that \ref{it1}) and \ref{it2}) force the set $\F$ to be
a factorization, it is easy to show that (\ref{log}) is a Lie
series whose valuation is strictly greater than 1. Hence, if $C$
is a conjugacy class of $\M(A,\theta)$, one has
\begin{equation}\label{lyndon}
\left(\underline C,\sum\log\underline{\langle
Y_i\rangle}\right)=\left(\underline C,\sum_{l\in Ly(A,\theta)
}\log{1\over1-l}\right)
\end{equation}
where $Ly(A,\theta)$ denotes the set of Lyndon traces (which is a
complete factorization of $\M(A,\theta)$ for the standard order
\cite{Lal}) and $(\ ,\ )$ is the scalar product  for which the
monomials form an orthonormal family.\\ If $C$ is not connected
(\ref{lyndon}) implies
\begin{equation}
\left(\underline
C,\sum_i\sum_m\frac1m\underline{Y_i^m}\right)=\left(\underline
C,\sum_l\sum_m\frac1ml^m\right)=0\nonumber
\end{equation}
and, the series  $\sum_i\sum_m\frac1m\underline{Y_i^m}$ being
positive, one gets $\langle Y_i\rangle\cap C=\emptyset$.\\ If $C$
is strongly primitive ({\it i.e.} $C$ is connected and $\sqrt
C=C$), equality (\ref{lyndon}) implies
\begin{equation}
\left(\underline
C,\sum_i\sum_m\frac1m\underline{Y_i^m}\right)=1.\nonumber
\end{equation}
Let $i$ and $m$ be such that $C\cap Y_i^m\neq\emptyset$. The
strong primitivity of $\underline C$ implies
\begin{equation}
\left(\underline C,\sum_i\sum_m\frac1m\underline{Y_i^m}\right)\geq
1\nonumber
\end{equation}
and the unicity of $Y_i$ follows. Furthermore, by $\Card C\cap
Y_i^m=m$, we show that $C\cap Y_i^m$ is a conjugacy class in
$\langle Y_i\rangle$.\\ Now, suppose  that $C$ is connected but
not primitive and let $p>1$ such that $C=\{g^p/g\in\sqrt C\}$. Let
$Y_i$ such that $\sqrt C\cap\langle Y_i\rangle\neq\emptyset$ is a
conjugacy class in $\langle Y_i\rangle$ (from the previous case,
$Y_i$ exists and is unique). One has $C\cap \langle
Y_i\rangle=C_j$ where each $C_j$ is a conjugacy class in $\langle
Y_i\rangle$. But for each $j$,
\begin{equation}
\left(\underline{C_j},\sum_m\frac1m\underline{Y_i}^m\right)=\frac1p\nonumber
\end{equation}
and from (\ref{lyndon}) one has
\begin{equation}
\left(\underline
C,\sum_m\frac1m\underline{Y_i}^m\right)=\frac1p.\nonumber
\end{equation}
The unicity of  $Y_i$ and $C_j$ follows.\\ \\ The proofs of the
two other implications are slight adaptations of the case where
$\theta=\emptyset$. The reader can refer to \cite{Sc} for more
details.\cqfd
%%%%%%%%%%%%%%%%%%%%%%%%%%%%%%%%%%%%%%%%%%%%%%%%
 Denote by ${\rm Cont}(\F)=\bigcup_{i\in J}Y_i$ the contents of the
  factorization $\F=(Y_i)_{i\in I}$. The following result
   is an extension of a classical result
due to Sch\"utzenberger \cite{onfact}.
\begin{corollary}
Let $\F$ be a complete factorization of $\M(A,\theta)$ (i.e. each
$Y_i$ is a singleton). Then for each conjugacy class $C$ we have
\[{\rm Card}(C\cap{\rm Cont}(\F))=\left\{\begin{array}{ll}1&\mbox{ if
} C \mbox{ is strongly connected }\\ &\mbox{(i.e. } C \mbox{ is
connected and } \sqrt C=C)\\ 0& otherwise
\end{array}\right.\]
\end{corollary}
Hence,  complete factorizations of $\M(A,\theta)$ receive the same
combinatorics than in the free case. In particular, one recovers
that the generating series of a complete factorization is equal to
the Hilbert series of the free partially commutative Lie algebra.
\section{Transitive Lazard sets\label{TLS}}

\subsection{Complete elimination strings and transitive Lazard sets}

Let ${\cal A}(2,A)$ be the free (non associative) algebra on $A$
whose product will be denoted by $(.,.)$ ({\it i. e.} the algebra of
the binary trees with leaves in $A$) .
\\
The canonical morphism ${\cal A}(2,A)\rightarrow \M(A,\theta)$
which is the identity on $A$ will be called the foliage morphism and
denoted by $f$ as in \cite{Re}. Remark that
$\theta_\M=\{(w,w')|ww'=w'w\mbox{ and
}\Alph(w)\cap\Alph(w')=\emptyset\}$ is a commutation relation on
$\M(A,\theta)$ (\cite{DL}).
\\
 In this section, we consider  a commutation alphabet
$(A,\theta)$ and some graphs whose vertices belong in ${\cal
A}(2,A)$ and such that the edges verify the property
 that $(t_1,t_2)$ is an edge if and only if $(f(t_1),
f(t_2))\in \theta_\M$.\\ Let $G=(V,E)$ be such a graph, we call a
{\bf Elimination String} ({\bf ES})  in $G$ a n-uplet of vertices
$(a_1,\cdots, a_n)$ such that for each $i\in[1,n ]$ and
$v_1,v_2\in V-\{a_1,\cdots,a_i\}$
\begin{equation}
(v_1,v_2)\in E\Rightarrow (v_1,a_i)\in E\mbox{ or } (v_2,a_i)\in E
\end{equation}
The vertex $a_1$ will be called the {\bf starting point} of the
{\bf ES}.
\\
An {\bf ES} $(a_1,\cdots,a_n)$ will be called {\bf complete} ({\bf
CES} in the following) if $E=\{a_1,\cdots, a_n\}$. A graph
admitting a {\bf CES} will be called {\bf type-H graph}.\\ Let
$n\geq 1$. In the sequel we denote ${\cal A}(2,A)^{\leq n}$ the
set of the trees with less than $n$ leaves.
\begin{example}
\begin{enumerate}
\item Let us consider the following graph
\[(A,\theta)=\begin{array}{ccccc}
a&-&d&-&e\\ |&&|&&\\ b&-&c&&
\end{array}\]
the family  $(a,d,b,c,e)$ is a {\bf CES} of $(A,\theta)$ and then
it is a type H graph.
\item The graph
\[
(A,\theta)=\begin{array}{ccc} a&-&b\\ c&-&d
\end{array}
\]
is not a type-H graph.
\end{enumerate}
\end{example}

\subsection{Transitive Lazard sets}
Let $v$ be a vertex of $G$, the {\bf H-star} of $G$ for $v$ to the
rank $n>0$ is the graph $G_n^{*v}=(V_n^{*v},E_n^{*v})$ defined by
\begin{equation}
\left\{
\begin{array}{l}
V_n^{*v}=(V\cap{\cal A}(2,A)^{\leq n}-v)\cup \{(v'v^m)\in{\cal
A}(2,A)^{\leq n}|m>0, (v',v)\not\in E\}\\ E_n^{*v}=\{(v_1,v_2)\in
(V_n^*)^2|(f(v_1),f(v_2))\in\theta_\M\}
\end{array}
\right. .
\end{equation}
\begin{example}
 For the following graph
\[(A,\theta)=a-b-c-d\]
then
\[(A^{*c}_4,\theta_4^{*c})=
\begin{array}{ccccccc}
&&ac&&\quad&\\ &&|&&&\\ a&-&b&-\quad ((a,c),c)&d\\ 
&&|&&&\\
&&(((a,c),c),c)&&&
\end{array}\]
\end{example}
\begin{proposition}\label{typeH}
Let $G$ be a type-H graph and $a_1$ be the starting point of a
$CES$ then $G_n^{*a_1}$ is a type-H graph.
\end{proposition}
{\bf Proof} Let $\lambda=(a_1,\cdots, a_m)$ be a {\bf CES} of $G$.
We construct a list $\Lambda$ of vertices of $G_n^{*a_1}=(V^*,E^*)$
removing $a_1$ of $\lambda$ and substituting each $a_i$ such that
$(a_i,a_1)\not\in E$ for the sequence
$$(a_i,a_1^{m-1}),\dots,(a_i,a_1),a_i.$$ Let us denote
$\Lambda=(a'_1,\cdots,a'_M)$ and suppose that $\Lambda$ is not a {\bf
CES} ; then it exists $\alpha<\beta<\gamma$ such that
$(a'_\beta,a'_\gamma)\in E^{*}$, $(a'_\alpha,a'_\beta)\not\in E^*$
and $(a'_\alpha,a'_\gamma)\not\in E^*$. Set
$a'_\alpha=(a_i,a_1^p)$, $a'_\beta=(a_j,a_1^q)$ and
$a'_\gamma=(a_k,a_1^r)$. The construction implies that $i\leq j<k$ and
$p=0$ or $q=0$. If $i=j$ then $p\neq 0$  which implies
$(a_i,a_1)\not\in E$ and $r=0$. Hence, $(a_1,a_k)\not\in E$ and
$q=0$. This implies that $(a_i,a_k)\in E$, $(a_1,a_i)\not\in E$
and $(a_1,a_k)\not\in E$ and contradicts the fact that $\lambda$
is a {\bf CES}. If $i< j$, suppose that $r=0$ (the case $q=0$ is
symmetric), we need to examine several cases
\begin{enumerate}
\item If $p=0$, then $q\neq 0$ (otherwise $(a_i,a_j),
(a_i,a_k)\not\in E$ and  $(a_j,a_k)\in E$ and this contradicts the fact that
$\lambda$ is a {\bf CES}). But, as $(a_k,a_i)\not\in E$,
$(a_j,a_k)\in E$ and as $\lambda$ is a {\bf CES}, one has
$(a_i,a_j)\in E$ and hence $(a_1,a_i)\not\in E$. Finally
$(a_1,a_j)\not\in E$  contradicts the fact that $\lambda$ is a
{\bf CES}.
\item If $p\neq 0$ and $q=0$, then, as $\lambda$ is a {\bf CES},
either $(a_i,a_j)\in E$ or $(a_i,a_k)\in E$. Suppose that
$(a_i,a_j)\in E$ (the other case is symmetric), one has
$(a_1,a_j)\in E$ (otherwise $(a'_\alpha,a'_\beta)\in E^*$). But,
$\lambda$ being a {\bf CES}, one gets $(a_1,a_k)\in E$ and by
$(a'_\alpha,a'_\gamma)\not\in E^*$, one obtains $(a_i,a_k)\not\in
E$. Finally $(a_1,a_j), (a_1,a_i)\not\in E$ and $(a_i,a_j)\in E$
contradicts the fact that $\lambda$ is a {\bf CES}.
\item If $p, q\neq 0$ then, as $(a'_\gamma,a'_\beta)\in E^*$, one
as $(a_1,a_k)\in E$ and by $(a'_\alpha,a'_\gamma)\not\in E^*$, one
obtains $(a_i,a_k)\not\in E$. Hence, $\lambda$ being a {\bf CES},
$(a_i,a_j)\in E$.  Hence $(a_1,a_j), (a_1,a_i)\not\in E$ and
$(a_i,a_j)\in E$ contradicts the fact that $\lambda$ is a {\bf
CES}.
\end{enumerate}
This prove that $\Lambda$ is a {\bf CES}.
 \cqfd
 Let $G=(A,\theta)$ be a
commutation alphabet considered as a graph and $L\subset{\cal
A}(2,A)$, we will say that $L$ is a {\bf transitive Lazard set} if
and only if for each $n> 0$, $L\cap{\cal A}(2,A)^{\leq
n}=\{s_1,\cdots, s_k\}$ such that it exists $k+1$ graphs
$G_1=(A_1,\theta_1),\dots,G_{k+1}=(A_{k+1},\theta_{k+1})$
satisfying the following conditions:
\begin{enumerate}
\item The first graph $G_1$ is equal to $G$.
\item The last graph $G_{k+1}$ is empty (i.e.
$G_{k+1}=(\emptyset,\emptyset)$)
\item For each $i<k+1$, $s_i\in A_i$ and $s_i$ is the starting
point of a {\bf CES} of $G_i$.
\item We have $G_{i+1}=(G_i)_n^{*s_i}$.
\end{enumerate}

\subsection{Type-H graphs and the transitive factorization theorem\label{TypeH}}

The classical properties of the (non commutative) Lazard sets hold
true and can be seen as consequences of the Transitive
Factorization theorem \cite{DL}. We recall them here.\\ 
Let
$(A,\theta)$ be a partially commutative alphabet and $B\subset A$.
Then $\M(B,\theta_B)$ is the left (resp. right) factor of a
bisection of $\M(A,\theta)$. Explicitly,
\[\M(A,\theta)=\M(B,\theta_B).\langle \beta_Z(B)\rangle\]
where $\langle \beta_Z(B)\rangle$ denotes the submonoid generated by
the set
\[\beta_Z(B)=\{zw/z\in Z, w\in \M(B,\theta_B), IA(zw)= \{z\}\}\]
and $IA(t)=\{z\in A|t=zw\}$ is the initial alphabet of the trace
$t$.\\ Let $B\subset A$, we say that $B$ is a {\bf transitively
factorizing subalphabet} (TFSA) if and only $\beta_Z(B)$ is a
partially commutative code \cite{DK2}. We proved the following theorem.
\begin{theorem}Duchamp-Luque\cite{DL}\label{transitive}
\begin{enumerate}
\item Let $B\subset A$. The following assertions are equivalent.
\begin{enumerate}
\item[(i)]The subalphabet $B$ is a TFSA.
\item[(ii)]The subalphabet $B$ satisfies the following condition.\\
For each $z_1\neq z_2\in Z$ and $w_1, w_2, w'_1,
w'_2\in\M(A,\theta)$ such that $IA(z_1w_1)=IA(z_1w'_1)=\{z_1\}$
and $IA(z_1w_2)=IA(z_2w'_2)=\{z_2\}$ we have
\[z_1w_1z_2w_2=z_2w'_2z_1w'_1\Rightarrow w_1=w'_1, w_2=w'_2. \]
\item[(iii)] For each $(z,z')\in Z^2\cap \theta$, the
dependence graph  ({\it ie} non-commutation) has no partial graph
like
\[z-b_1-\dots-b_n-z'.\]
with $b_1,\dots,b_n\in B$.
\end{enumerate}
\item Let $(B,Z)$ be a partition of $A$
\begin{quote}
(i) We have the decomposition
\[L(A,\theta)=L(B,\theta_B)\oplus J\]
where $J$ is the Lie ideal generated (as a Lie algebra) by
\[\tau_Z(B)=\{[ \dots [z,b_1], \dots b_n]\quad|\quad zb_1 \dots b_n\in
\beta_Z(B)\}. \] (ii) The subalgebra $J$ is a free partially
commutative Lie algebra if $B$ is a TFSA of $A$.\\ (iii)
Conversely if $J$ is a free partially commutative Lie algebra with
code $\tau_Z(B)$ then %%@
$B$ is a TFSA.
\end{quote}
\end{enumerate}
\end{theorem}
Applying theorem \ref{transitive} and taking the inductive limit
of the process one obtains.
\begin{proposition}
Let $L$ be a transitive Lazard set.
\begin{enumerate}
\item The foliage $f(L)$ is a complete factorization of
$\M(A,\theta)$.
\item Let $\Pi$ be the unique morphism
${\cal A}(2,A)\rightarrow L(A,\theta)$ such that $\Pi(a)=a$ for
each letter $a\in A$. Then $\Pi(L)$ is a basis of the Free Lie
algebra $L(A,\theta)$.
\end{enumerate}
\end{proposition}
Such a factorization will be called Transitive Lazard
Factorization (TLF).  Not all the trace monoids possess a TLF. For
example, in the graph $a-b\atop c-d $ we can not find a {\bf CES}.
Nevertheless, the property ``having a TLF" is decidable as shown
by the following result.
\begin{theorem}
A trace monoid admits a TLF if and only if its commutation
alphabet is a type-H graph.
\end{theorem}
{\bf Proof} It suffices to remark that a trace monoid has a TLF if
and only if on can construct a transitive Lazard set from its
commutation graph which is a consequence of proposition
\ref{typeH}.\cqfd

\section{Transitive Hall sets}

Let us define a {\bf Transitive Hall Set} ({\bf THS}) $H$ as a
family of trees $(h)_{h\in H}$ endowed with a total order $<$ such
that
\begin{enumerate}
\item The family $(f(h))_{h\in H}$ is a complete factorization of
$\M(A,\theta)$ (for the reverse order).
\item The set $H$ contains the alphabet $A$.
\item If $h=(h',h'')\in H-A$ then $h''\in H$ and $h<h''$.
\item If $h=(h',h'')\in{\cal A}(2,A)-A$ then $h\in H$ if and only if the
four following assertions are true :
\begin{enumerate}
\item The two sub-trees $h'$ and $h''$ belong to $H$.
\item We have the inequality $h'<h''$.
\item The foliage of the two sub-trees are not related by the ``disjoint commutation'' relation (i.e. $(f(h'),f(h''))\not\in\theta_\M$).
\item Either $h'\in A$ or $h'=(x,y)$ with $y\geq h''$.
\end{enumerate}
\end{enumerate}
Note that in the classic non-commutative theory of Hall set, one gets (1)
from the axioms (2), (3) and (4). Here, it is not the case. For
example, if we consider the commutation graph
\[a\qquad b-c\]
and a set $H$ such that
\[\begin{array}{rcl}
H\cap{\cal A}(2,A)^{\leq 3}
&=&\{a,(b,a),((c,a),a),c,((c,a),c),(c,a),(b,(c,a)),\\ &&
  ((b,a),c),b,((b,a),b),(b,a)\}.
\end{array}\]
We suppose that  trees above are ordered from the right to the
left. Clearly,
 trees listed here verify  axioms (2), (3) and (4) of the definition of
 transitive Hall sets, but we can not complete $H\cap {\cal A}(2,A)^{\leq 3}$
to construct a complete factorization as $c.ba$ is a decreasing
decomposition of the trace $bca$ which belongs to $f(H)$.

As in the free case we have a perfect correspondence between the
notions of transitive Hall and Lazard sets.
\begin{theorem}
A set $L$ is a transitive Lazard set if and only if it is a
transitive Hall set.
\end{theorem} {\bf Sketch of the proof} A slight adaptation of the non commutative case (see \cite{Re})
shows that each transitive Lazard set is a Transitive Hall set.\\
Let us sketch a proof of the converse. The notion of Transitive Lazard sets
is independent of the filtration in the following sense: $L$ is a
Lazard set if and only if for each finite closed set $E$ ({\it ie}
$(t_1,t_2)\in E\Rightarrow t_1,t_2\in E$) denoting $L\cap
E=\{s_1,\dots,s_k\}$ there exist $k+1$ graphs $G_i=(A_i,\theta_i)$
verifying
\begin{enumerate}
\item $G_1=(A,\theta)$,
\item $G_{k+1}=(\emptyset,\emptyset)$
\item For each $i<k+1$, $s_i\in A_i$ and $s_i$ is the starting
point of a CES over $G_i$.
\item $A_{i+1}=(A_i-s_i)\cup \{(bc_i^n)\in E|n>0,b\in A_i \mbox{ and }0 (b,c_i)\not\in\theta_i\}$
and $\theta_{i+1}=\theta_{A_{i+1}}$.
\end{enumerate}
The fact that, for each closed set $E$,  any transitive Hall set
satisfies the four previous properties will be proved by induction on the cardinal
of $E$. Let $H$ be a THS. If ${\rm Card\  }E=1$ the result is
obvious. We suppose now that ${\rm Card\  }E>1$ and we set
\begin{enumerate}
\item[] $c=max\{h\in H\cap E\}$
\item[] $X=\{(ac^n)|a\in A-c,n\geq
0,(a,c)\in\theta\}\cap\{b|(b,c)\in\theta\}$
\item[] $\theta_X=\theta_\M\cap X\times X$
\item[] $H'=H\cap H\cap {\cal A}(2,X)$ and $E'=E\cap {\cal A}(2,X)$
\end{enumerate}
As in \cite{Re}, $H'$ is a  transitive Hall set for the alphabet
$X$ and  $H'\cap E'= H\cap E'$. By induction $H'$ is a transitive
Lazard set and one  constructs $k+1$ graphs
$G_1=(A_1,\theta_1),\dots, G_{k+1}=(A_{k+1},\theta_{k+1})$
verifying (1), (2), (3), and (4). Remark that if $H\cap
E'=\{s_1,\dots,s_n\}$, one has $H\cap E=\{s_0=c,s_1,\dots,s_n\}$.
 Suppose that it exists $(a,b)\in\theta$ such that $(a,c)\not\in
\theta$ and $(b,c)\not\in\theta$ then a quick examination shows
that $(f(h))_\H$ is not a factorization. Hence, $c$ is the
starting point of a CES in $G_0=(A,\theta)$ and we have
$A_1=(A-c)\cup \{(ac^n)\in E|n>0,a\in A \mbox{ and }0
(a,c)\not\in\theta\}$. It follows that $H$ and $E$ verify the
assertions (1), (2), (3) and (4). This proves that $H$ is a
transitive Lazard set.
  \cqfd
 This correspondence is very useful to
construct decomposition algorithms. We can construct {\bf standard
sequences} of Hall trees $(h_1,\cdots,h_n)$ such that
$L(A,\theta)$ for each $i\in[1,n]$ either $h_i\in A$, or
$h_i=(h'_i,h''_i)$ with $h''_i\geq h_{i+1},\dots,h_n$. In a
standard sequence an {\bf ascent} is an index $i$ such that
$h_i<h_{i+1}$ and a {\bf legal  ascent} is an ascent $i$ such that
$h_{i+1}\geq h_{i+2},\cdots, h_n$ (these definitions are due to
Sch\"utzenberger \cite{Sc}).  Let $s$ be a standard sequence and
 $i$ a  legal
ascent. We write $s\rightarrow s'$ if
$s'=(h_1,\cdots,h_{i-1},(h_i,h_{i+1}),h_{i+2},\cdots, h_n)$ when
$(f(h_i),f(h_{i+1}))\not\in\theta_\M$ and
$s'=(h_1,\cdots,h_{i-1},h_{i+1},h_i,h_{i+2},\cdots,h_n)$
otherwise. The transitive closure
$\displaystyle\mathop\rightarrow^*$ of $\rightarrow$ is such that, 
for each standard sequence it exists a unique decreasing standard
sequence $s'$ such that $s \displaystyle\mathop\rightarrow^*s'$.
Using this property on a sequence of letters, we obtain an
algorithm which allows to find the factorization of a trace in a
decreasing concatenation of Hall traces.
\begin{example}\label{HallEx2}
Let us consider the following commutation alphabet:
\[(A,\theta)=a-b-c-d .\]
and let $H$ be a transitive hall set such that\\ \\ $H\cap{\cal
A}(2,A)^{\leq
3}=\{c,b,a,(a,c),((a,c),c),(a,(a,c)),(d,b),((d,b),b),(d,(b,a)),$\\
$(d,(a,c)),(d,a),((d,a),a),(d,(d,b)),(d,(d,a)),d\}.$\\ \\ We can
compute the factorization of the word $bcaccbdbddad$ in the
following way $$
\begin{array}{c}
(b,c,a,c,c,b,d,b,d,d,a,d)\\ \downarrow\\
(b,c,a,c,c,b,d,b,d,(d,a),d)\\ \downarrow\\
(b,c,a,c,c,b,d,b,(d,(d,a)),d)\\ \downarrow\\
(b,c,(a,c),c,b,(d,b),(d,(d,a)),d)\\ \downarrow\\
(b,c,((a,c),c),b,(d,b),(d,(d,a)),d)\\ \downarrow\\
(b,c,b,((a,c),c),(d,b),(d,(d,a)),d)\\ \downarrow\\
(c,b,b,((a,c),c),(d,b),(d,(d,a)),d)
\end{array}
$$ which gives $bcaccbdbddad=c.b.b.acc.db.dda.d$.
\end{example}
Let $s=(h_1,\cdots,h_n)$ be a standard sequence and $i$ a legal
ascent, we define
\begin{equation}
\lambda_i(s)=(h_1,\dots,h_{i-1},(h_i,h_{i+1}),h_{i+2},\dots,h_n)
\end{equation}
and
\begin{equation}
\rho_i(s)=(h_1,\cdots,h_{i-1},h_{i+1},h_i,h_{i+2},\cdots,h_n).
\end{equation}
The {\bf derivation tree} of $s$ is the tree $T(s)$ satisfying the
following
\begin{enumerate}
\item if $s$ is a decreasing sequence then $T(s)$ is only the root
labeled $s$
\item otherwise, we consider the greatest legal ascent $i$ of $s$.
Then
\begin{enumerate}
\item if $(h_i,h_{i+1})\in \theta_\M$, the root of the tree $T(s)$ is $s$
and $T(s)$ has only one sub-tree $T(\rho_i(s))$.
\item otherwise, the root of $T(s)$ is $s$, the left sub-tree of
$T(s)$ is $T(\lambda_i(s))$ and the right sub-tree of $T(s)$ is
$T(\rho_i(s))$
\end{enumerate}
\end{enumerate}
If we denote $[s]=[h_1]\cdots [h_n]$, one obtains
\begin{equation}
[s]=\sum_{s'\in{\goth F}(T(s))}[s']
\end{equation}
where ${\goth F}(T(s))$ denotes the set of the leaves of $T(s)$.
Applying this equality to sequences of words, one gets an
algorithm allowing to decompose a polynomial in the PBW  basis
associated to a transitive Hall set.
\begin{example} We use the transitive Hall set defined in the example \ref{HallEx2}. One
has for example:
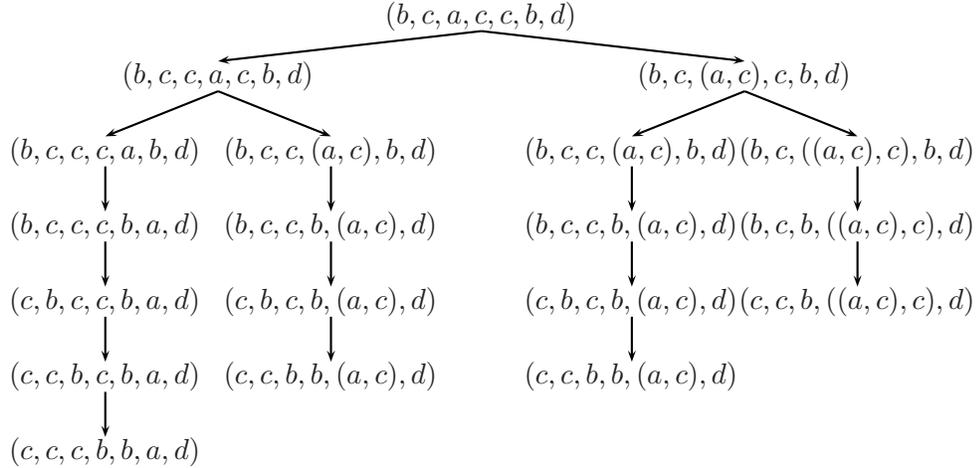
\begin{center}\small
\begin{pspicture}(13,6)
%\psgrid
\rput(6.5,5.8){$(b,c,a,c,c,b,d)$} \rput(3,5){$(b,c,c,a,c,b,d)$}
\rput(1.5,4){$(b,c,c,c,a,b,d)$} \rput(1.5,3){$(b,c,c,c,b,a,d)$}
\rput(1.5,2){$(c,b,c,c,b,a,d)$} \rput(1.5,1){$(c,c,b,c,b,a,d)$}
\rput(1.5,0){$(c,c,c,b,b,a,d)$} \rput(4.5,4){$(b,c,c,(a,c),b,d)$}
\rput(4.5,3){$(b,c,c,b,(a,c),d)$}
\rput(4.5,2){$(c,b,c,b,(a,c),d)$}
\rput(4.5,1){$(c,c,b,b,(a,c),d)$} \rput(10,5){$(b,c,(a,c),c,b,d)$}
\rput(8.5,4){$(b,c,c,(a,c),b,d)$}
\rput(8.5,3){$(b,c,c,b,(a,c),d)$}
\rput(8.5,2){$(c,b,c,b,(a,c),d)$}
\rput(8.5,1){$(c,c,b,b,(a,c),d)$}
\rput(11.5,4){$(b,c,((a,c),c),b,d)$}
\rput(11.5,3){$(b,c,b,((a,c),c),d)$}
\rput(11.5,2){$(c,c,b,((a,c),c),d)$} \psline{->}(6.5,5.6)(3,5.2)
\psline{->}(3,4.8)(1.5,4.2) \psline{->}(1.5,3.8)(1.5,3.2)
\psline{->}(1.5,2.8)(1.5,2.2) \psline{->}(1.5,1.8)(1.5,1.2)
\psline{->}(1.5,0.8)(1.5,0.2) \psline{->}(3,4.8)(4.5,4.2)
\psline{->}(4.5,3.8)(4.5,3.2) \psline{->}(4.5,2.8)(4.5,2.2)
\psline{->}(4.5,1.8)(4.5,1.2) \psline{->}(6.5,5.6)(10,5.2)
\psline{->}(10,4.8)(8.5,4.2) \psline{->}(8.5,3.8)(8.5,3.2)
\psline{->}(8.5,2.8)(8.5,2.2) \psline{->}(8.5,1.8)(8.5,1.2)
\psline{->}(10,4.8)(11.5,4.2) \psline{->}(11.5,3.8)(11.5,3.2)
\psline{->}(11.5,2.8)(11.5,2.2)
\end{pspicture}\end{center}
Which allows to write 
$$
bcaccbd=c.c.c.b.b.a.d+2c.c.b.b.[a,c].d+c.b.b.[[a,c],c].d .
$$

\end{example}
\section{Conclusion}

The stable concept of transitive factorization allows the
adaptation of the Hall machinery as it   has been explicited here.
The construction is characteristic free. It would be interesting
to investigate such constructions in the case of $p$-commutations
\cite{DDM,DL2}.

\end{document}